# Rejoinder

**Javier M. Moguerza and Alberto Muñoz**

## 1. INTRODUCTION

We are very grateful to the Executive Editors George Casella and Edward George for their active interest in our paper and for organizing this challenging discussion. We also thank all the discussants for their insightful and stimulating comments.

When we submitted the original manuscript in 2003, we were tempted to go for a more general paper on kernel methods. We decided to focus on support vector machines (SVMs), waiting for a mature development of the new and exciting ideas related to kernel methods, such as manifold learning and other related topics. We will refer to some of these methods below. Let us begin, first, with some general considerations.

Regarding the question of the dimensionality induced by the feature space, Hastie and Zhu remark in their comment that usual kernels do not automatically lead to infinite-dimensional feature spaces. They give a nice example that involves the radial (Gaussian) kernel function. This agrees with results in Keerthi and Lin [8], where an explanation of the performance of the Gaussian kernel is given when, according to the notation in the comment by Hastie and Zhu, $\gamma \to 0$ and $\lambda$ is chosen in the appropriate way. In this case, the SVM classifier converges to a linear SVM classifier, and the effective dimension of the kernel is finite, agreeing with the empirical conclusion provided by the discussants.

We also agree with the assertions of some of the discussants regarding the probabilistic interpretability of the SVM output (the sign of some estimated function). Our comment was rather along the line of Sollich [18], who proposed to make Bayesian methods available for the support vector methodology, while leaving as much as possible of the standard SVM framework intact. This is not an easy task. In fact, as Bartlett, Jordan and McAuliffe remark, sparseness and the precise estimation of conditional probabilities are hard to reconcile.

Regarding the role of differentiability in SVMs (misplaced in the opinion of Bartlett, Jordan and McAuliffe), it is convenient to recall that the differentiable formulation of the SVM problem allows its solution by the use of standard Newton-type methods for convex optimization. Under the availability of second order derivatives (and this is the case for SVMs), these methods are known to be the most efficient ones for the solution of smooth problems.

We thank some of the discussants for turning the attention of the reader to general kernel methods. In particular, we appreciate the Bartlett, Jordan and McAuliffe effort to make clearer the potential impact of reproducing kernel Hilbert space (RKHS) methods. Regarding the origins of RKHS in statistics, for the sake of completeness, we strongly recommend reading the conversation with Emanuel Parzen in [14].

Given the history of SVMs, perfectly outlined by Wahba in the introduction of her comment, we do not like to think of SVMs as a "modest" variant of some standard statistical methodology (as suggested by Bartlett, Jordan and McAuliffe). Using a similar (a posteriori) reasoning, some strict mathematicians might think that RKHS methods in statistics are just a small variation on the general theory of Hilbert spaces. Of course, this is far from true. We rather think that the support vector methodology, followed closely by kernel methods, has been able to synthesize a variety of techniques from different fields, leading to a more unified framework for learning theory [5]. In addition, the geometrical viewpoint of SVMs allows new approaches to long-familiar problems, as illustrated in the next section.

## 2. KERNEL METHODS REVISITED

One interesting point regarding the geometrical interpretation of SVMs is that they have stirred the development of new techniques driven by the geometrical properties of the kernel. Some of these techniques have not so far been mentioned in the discussion. We now briefly describe two relevant examples.

### 2.1 One-Class SVMs

An example of a new method that has arisen from a geometrical point of view is one-class SVMs [16]. One-class SVMs deal with a problem related to estimating high density regions from data samples. The method computes a binary function that takes the value $+1$ in "small" regions that contain most data points and takes the value $-1$ elsewhere. The strategy of the one-class support vector method is to map the data points into the feature space determined by a kernel function and to calculate a hyperplane that separates the mapped data $\{\Phi(\mathbf{x}_i)\}_{i=1}^n$ from the origin, where $\Phi$ is the mapping induced by the kernel function. With this aim, the one-class SVM algorithm solves the quadratic optimization problem

$$(2.1) \quad \begin{aligned} \min_{\mathbf{w},b,\xi} \quad & \frac{1}{2}\|\mathbf{w}\|^2 - b + \frac{1}{\nu n}\sum_{i=1}^n \xi_i \\ \text{s.t.} \quad & \mathbf{w}^T\Phi(\mathbf{x}_i) \geq b - \xi_i, \quad i=1,\ldots,n, \\ & \xi_i \geq 0, \quad i=1,\ldots,n, \end{aligned}$$

where $\xi_i$ are slack variables, $\nu \in [0,1]$ is an a priori fixed constant which represents the fraction of outlying points and $b$ is the decision value which determines whether a given point belongs to the estimated high density region. The decision function will take the form $h(\mathbf{x}) = \text{sign}(\mathbf{w}^{*T}\Phi(\mathbf{x}) - b^*)$, where $\mathbf{w}^*$ and $b^*$ are the values of $\mathbf{w}$ and $b$ at the solution of problem (2.1). The hyperplane $\mathbf{w}^{*T}\Phi(\mathbf{x}) - b^* = 0$ separates from the origin the mapped data for which the decision function $h(\mathbf{x}) = +1$. Problem (2.1) is smooth and convex, and follows the SVM idea of building a hyperplane in a feature space.

It is apparent that solving the problem of estimating high density regions by building a separating hyperplane in a feature space is not trivial. Next, we provide an original statistical explanation of one-class SVMs. Consider the class of real-valued functions

$$(2.2) \quad \mathcal{G} = \{g > 0 \,|\, \forall \mathbf{x}, \mathbf{y} \in X, \\ g(\mathbf{x}) > g(\mathbf{y}) \iff f(\mathbf{x}) > f(\mathbf{y})\},$$

where $X$ is the input space and $f$ is the data density function. To estimate the outlying points that correspond to the proportion $\nu$, all we have to do is use the order induced by any function $g \in \mathcal{G}$ on the data sample $\{\mathbf{x}_1,\ldots,\mathbf{x}_n\}$. This is equivalent to solving the optimization problem

$$(2.3) \quad \begin{aligned} \max_{\lambda} \quad & -\sum_{i=1}^n \lambda_i g(\mathbf{x}_i) \\ \text{s.t.} \quad & \sum_{i=1}^n \lambda_i = 1, \\ & 0 \leq \lambda_i \leq \frac{1}{\nu n}, \quad i=1,\ldots,n, \end{aligned}$$

where, at the solution $\lambda^* = (\lambda_1^*,\ldots,\lambda_n^*)^\mathbf{T}$, $\lambda_i^* > 0$ if $\mathbf{x}_i$ is an outlying point [i.e., $\lambda_i^* > 0$ for small values of $g(\mathbf{x}_i)$] and $\lambda_i^* = 0$ otherwise. The dual of this linear problem is

$$(2.4) \quad \begin{aligned} \min_{b,\xi} \quad & -b + \frac{1}{\nu n}\sum_{i=1}^n \xi_i \\ \text{s.t.} \quad & g(\mathbf{x}_i) \geq b - \xi_i, \quad i=1,\ldots,n, \\ & \xi_i \geq 0, \quad i=1,\ldots,n. \end{aligned}$$

It can be shown that, at the solution, the values of the objective functions of problems (2.3) and (2.4) coincide (see, e.g., [1]). Moreover, the solution of problem (2.3) can be straightforwardly calculated from the solution of problem (2.4).

The one-class SVM problem (2.1) and problem (2.4) are very similar. It becomes clear now that the solution of problem (2.1) (the one-class SVM solution) tries to estimate a function $g \in \mathcal{G}$ by $\hat{g}(\mathbf{x}) = \mathbf{w}^{*T}\Phi(\mathbf{x})$, that is, by estimating in the feature space the weights $\mathbf{w}^*$ of a hyperplane with minimum norm. This is achieved through the inclusion of the term $1/2\|w\|^2$ in the objective function of problem (2.1).

Appropriate mappings and kernels to solve the problem of estimating high density regions (density level sets) using one-class SVMs are derived and can be consulted in [13].

### 2.2 Combination of Kernels

Another example of a technique developed from geometric considerations of the kernel is now illustrated. This method falls in the category of "further advances" mentioned by Bousquet and Schölkopf at the end of their comment. In particular, we build, for classification purposes, what they call a joint kernel (mixing inputs and outputs). This joint kernel is built by the combination of a set of kernels. A key point of our proposal is that the constructed kernel tries to capture the "right" notion of similarity. This agrees with the comment by Bousquet and Schölkopf about the relationship between the good performance of SVMs in practice and the appropriate prior knowledge about the problems incorporated by kernels. Thus, we will work with similarity matrices instead of kernel matrices. In fact, as Wahba points out, Euclidean distances (and therefore similarities) can be derived from positive definite kernels.

The idea, in geometric terms, is introduced next. If kernels are being used, points in a sufficiently small neighborhood in the feature space should belong to the same class (excluding points very close to the decision surface). As a consequence, if we are going to classify a data set by relying on a given similarity matrix, points close to each other using such similarities should, in general, be in the same class. Therefore, we have to construct a similarity matrix $K^*$ with entries $K^*(x_i, x_j)$ that are large for $x_i$ and $x_j$ in the same class (i.e., $y_i = y_j$), and small for $x_i$ and $x_j$ in different classes. For instance, if two kernels $K_1$ and $K_2$ are to be combined, a possible choice is

$$(2.5) \quad K^*(x_i, x_j) = \begin{cases} \max(K_1(x_i, x_j), K_2(x_i, x_j)), \\ \quad \text{if } y_i = y_j, \\ \min(K_1(x_i, x_j), K_2(x_i, x_j)), \\ \quad \text{otherwise.} \end{cases}$$

It is immediate to show that (2.5) is equivalent to

$$(2.6) \quad K^* = \tfrac{1}{2}(K_1 + K_2) + \tfrac{1}{2}Y|K_1 - K_2|Y,$$

where $Y = \text{diag}(y)$ is a diagonal matrix whose nonzero elements are the data labels, that is, $y_i \in \{-1, +1\}$.

Let $K_1, K_2, \ldots, K_M$ be the available set of $M$ input kernel matrices, all of which are obtained from the same data sample $\{\mathbf{x}_1, \ldots, \mathbf{x}_n\}$. The extension of the previous idea to the combination of more than two kernel matrices is

$$(2.7) \quad K^* = \bar{K} + Y \sum_{i<j} g(K_i - K_j) Y,$$

where $\bar{K}$ is the average of the kernel matrices and $g$ is a function that quantifies the difference of information between kernel matrices. The function $g$ must have the property that if $K_i$ and $K_j$ tend to produce the same classification results, then $g(K_i - K_j)$ should be almost null. A particular case of the previous equation is

$$(2.8) \quad K^* = \bar{K} + Y \sum_{m=1}^{M} |K_m - \bar{K}| Y.$$

This and other choices for $g$ can be consulted in [12].

Next, we show how this method (denoted AV for absolute value) can be used to improve the performance of single kernels. With this aim, we will use the breast cancer data set, made up of 683 observations with 9 features each [11]. We have considered three kernels: a polynomial kernel $K_1(x, z) = (1 + x^T z)^2$, a Gaussian kernel $K_2(x, z) = \exp^{(-\|x-z\|^2)}$ and a linear kernel $K_3(x, z) = x^T z$. We will compare SVMs using these kernels with the AV combination method and a semidefinite programming (SDP) technique for building linear combinations of kernels developed by Lanckriet et al. [9]. The data set has been randomly partitioned ten times into a training set and a test set, and for each method, a run of the experiment has been done over each partition. The average results are shown in Table 1. The AV method provides the best results (a test error of 3.1%), using significantly less support vectors than the other methods. The SDP method improves only the results of the Gaussian and the polynomial kernel.

TABLE 1
*Percentage of misclassified data and support vectors for the cancer data. Standard deviations in parentheses*

|  | Training error | Test error | % SV |
|---|---|---|---|
| $K_1$: Polynomial | 0.1 (0.1) | 7.8 (2.5) | 8.3 (0.8) |
| $K_2$: Gaussian | 0.0 (0.0) | 10.8 (1.7) | 65.6 (1.0) |
| $K_2$: Linear | 2.6 (0.5) | 3.7 (1.8) | 7.1 (0.8) |
| AV | 2.4 (0.3) | 3.1 (1.3) | 2.9 (0.4) |
| SDP | 0.0 (0.0) | 6.2 (1.6) | 65.5 (1.9) |

2.2.1 *Parameter selection.* Techniques for the combination of kernels can be successfully applied to the problem of parameter selection in kernel methods. This links with the comment of Bousquet and Schölkopf about the need for further research on satisfactory alternatives, other than cross-validation, to choose the parameters in kernel methods. We illustrate this situation using a collection of Gaussian kernels on the cancer data set. Let $\{K_1, \ldots, K_{12}\}$ be a set of Gaussian kernels $K_c(\mathbf{x}, \mathbf{y}) = \exp^{(-\|\mathbf{x}-\mathbf{y}\|^2/c)}$ with parameters $c = 0.1, 1, 10, 20, 30, 40, 50, 60, 70, 80, 90$ and $100$, respectively. This wide set covers a realistic range of possible values for the kernel parameter. The test errors for 12 SVMs using these Gaussian kernels range from 3.1% to 24.7%. In this case, the AV method, combining the 12 Gaussian kernels, gives the best result obtained using only one of the Gaussian kernels under consideration (with a test error of 3.1%). It is important to note that the performance of the AV method (which is parameter-free) is not affected by the inclusion of kernels with a bad generalization performance. Since, in general, the best parameter choice is not known in advance, the methodology just described provides an alternative that minimizes the effect of bad parameter selection.

## 3. THE BIAS–VARIANCE PROBLEM

Regarding the comments on statistical consistency provided by Bartlett, Jordan and McAuliffe, also pointed out by Bousquet and Schölkopf, we agree that the Vapnik–Chervonenkis (VC) dimension is not central in the analysis of SVMs. In fact, in [6], the bias–variance problem is analyzed in the context of regularization for the quadratic loss function (the analysis by Steinwart cited by the discussants is for the $L_1$-SVM, as Bousquet and Schölkopf remark). Cucker and Smale [6] replaced the VC dimension by the radius $r$ of a ball in a RKHS space ($r$ is the norm in the RKHS of the minimizer of the empirical risk). Since the regularization parameter $\lambda$ (using the notation of Bartlett, Jordan and McAuliffe) is inversely proportional to $r$, large values of $\lambda$ correspond to large bias, while small values of $\lambda$ lead to large variance. The Cucker and Smale paper [6] also contains a theorem (Corollary 2) that is in agreement with the discussants' comment about the fact that the regularization coefficient must decrease with the sample size. In addition, statistical consistency can be derived from the results in the paper if a rich enough kernel is used (i.e., a universal kernel, in the sense of Steinwart).

## 4. DIFFERENTIAL GEOMETRY METHODS AND KERNEL METHODS

In her comment, Wahba introduces a particular method for learning the kernel data matrix for the purpose of manifold unfolding. As Wahba and her co-authors remark in [10], the manifold unfolding problem is closely related to the construction of a kernel. In fact, Ham, Lee, Mika and Schölkopf [7] showed that several of the proposed techniques for manifold learning, namely ISOMAP [19], graph Laplacian eigenmap [2] and locally linear embedding [15], can be interpreted as particular cases of kernel principal component analysis [17]. As Ham and co-authors point out, these techniques can be viewed as a "warping of the input space into a feature space where the manifold is flat."

In this regard, Burges [4] described the intrinsic geometry of the manifold which arises for a particular choice of the kernel. In particular, he shows that the Riemannian metric induced on the manifold by its embedding can be expressed, in terms of the kernel, in closed form. A closely related approach can be found in [20]. It is worth mentioning that the implicit geometric assumption within manifold unfolding is that the decision surface (for the case of classification) is smooth with respect to the underlying geometry [3].

Finally, we would like to thank David Rios and Francisco J. Prieto for their careful reading of the manuscript and suggestions. We hope that our paper and this interesting discussion encourage the statistical community to pursue further research on support vector machines and other related methodologies.